\documentclass[11pt]{article}

\date{21st October 2024}

\title{\sc Constructing non-AMNM weighted convolution algebras for every semilattice of infinite breadth}
\author{\sc Yemon Choi, Mahya Ghandehari, Hung Le Pham}

\usepackage[a4paper,left=30mm,right=30mm,top=30mm,bottom=30mm]{geometry}

\usepackage{hyperref}
\hypersetup{
colorlinks=true,
linkcolor=red, 
citecolor=blue,
}

\usepackage{sectsty}
\allsectionsfont{\sffamily}

\RequirePackage{amsmath,amssymb, graphicx}
\numberwithin{equation}{section}

\usepackage{amsthm}
\newcounter{pulse}[section]
\numberwithin{pulse}{section}

\newcommand{\newheadfont}{\sc} 

\newcommand{\thmsep}{\topsep} 

\newtheoremstyle{newplain} 
  {\thmsep}   
  {\thmsep}   
  {\itshape}  
  {0pt}       
  {\newheadfont} 
  {.}         
  {5pt plus 1pt minus 1pt} 
  {}          

\newtheoremstyle{newdef} 
  {\thmsep}   
  {\thmsep}   
  {\normalfont}  
  {0pt}       
  {\newheadfont} 
  {.}         
  {5pt plus 1pt minus 1pt} 
  {}          

\newtheoremstyle{newrem}
  {\thmsep}   
  {\thmsep}   
  {\normalfont}  
  {0pt}       
  {\newheadfont} 
  {.}         
  {5pt plus 1pt minus 1pt} 
  {}          

\theoremstyle{newplain}
\newtheorem{thm}[pulse]{Theorem}
\newtheorem{prop}[pulse]{Proposition}
\newtheorem{lem}[pulse]{Lemma}
\newtheorem{claim}[pulse]{Claim}

\theoremstyle{newdef}
\newtheorem{dfn}[pulse]{Definition}

\newtheorem{eg}[pulse]{Example}

\theoremstyle{newrem}
\newtheorem{rem}[pulse]{Remark}

\newcounter{alphacount}

\theoremstyle{newplain}
\newtheorem{thmalpha}[alphacount]{Theorem}

\newenvironment{newnum}{%
\begin{enumerate}

}{\end{enumerate}\ignorespacesafterend}

\newcommand{\veps}{\varepsilon}
\newcommand{\om}{\omega}
\newcommand{\lm}{\lambda}

\newcommand{\defeq}{\mathbin{:=}}
\newcommand{\dt}[1]{\textit{#1}\/}  

\newcommand{\blob}{\bullet}

\newcommand{\set}[1]{\left\{#1\right\}}
\newcommand{\abs}[1]{{\left\vert#1\right\vert}}

\newcommand{\Cplx}{{\mathbb C}}
\newcommand{\Nat}{{\mathbb N}}

\newcommand{\Mult}{\operatorname{Mult}}

\newcommand{\divides}{\!\mid\!}
\newcommand{\fac}{\operatorname{fac}} 
\newcommand{\FBP}{\operatorname{FBP}} 
\newcommand{\filgen}[1]{\operatorname{fil}(#1)}  

\newcommand{\ssfont}[1]{{\mathsf{#1}}}
\newcommand{\ssa}{\ssfont{a}}
\newcommand{\ssb}{\ssfont{b}}

\newcommand{\ssd}{\ssfont{d}}
\newcommand{\ssg}{\ssfont{g}}

\newcommand{\ssp}{\ssfont{p}}

\newcommand{\ssx}{\ssfont{x}}
\newcommand{\ssy}{\ssfont{y}}
\newcommand{\ssz}{\ssfont{z}}

\newcommand{\cA}{{\mathcal A}}
\newcommand{\cB}{{\mathcal B}}
\newcommand{\cC}{{\mathcal C}}

\newcommand{\cE}{{\mathcal E}}
\newcommand{\cF}{{\mathcal F}}
\newcommand{\cG}{{\mathcal G}}

\newcommand{\cP}{{\mathcal P}}

\newcommand{\cR}{{\mathcal R}}
\newcommand{\cS}{{\mathcal S}}
\newcommand{\cT}{{\mathcal T}}

\newcommand{\cX}{{\mathcal X}}

\newcommand{\FIN}{\operatorname{\mathcal{P}}\nolimits^{\rm fin}} 

\newcommand{\dotcup}{\mathbin{\dot{\cup}}}

\newcommand{\join}{\operatorname{join}}

\newcommand{\TMAX}{{\cT}_{\rm max}}
\newcommand{\TMIN}{{\cT}_{\rm min}}
\newcommand{\TORT}{{\cT}_{\rm ort}}
\newcommand{\level}{\operatorname{level}}

\newcommand{\sscap}{\wedge} 
\newcommand{\sscup}{\vee} 

\newcommand{\bigsscup}{\bigvee} 

\newcommand{\setcap}{\cap} 
\newcommand{\setcup}{\cup} 
\newcommand{\bigsetcap}{\bigcap} 
\newcommand{\bigsetcup}{\bigcup} 

\begin{document}
\maketitle

\hfil\textit{Dedicated to the memory of H. Garth Dales (1944--2022)}\hfil

\begin{abstract}
The AMNM property for commutative Banach algebras is a form of Ulam sta\-bili\-ty for multiplicative linear functionals. We show that on any semilattice of infinite breadth, one may construct a weight for which the resulting weighted convolution algebra fails to have the AMNM property.
Our work is the culmination of a trilogy started in \cite{CGP-SGF} and continued in \cite{CGP-EJC}. In particular, we obtain a refinement of the main result of \cite{CGP-EJC}, by establishing a dichotomy for union-closed set systems that has a Ramsey-theoretic flavour. 


\end{abstract}

\begin{section}{Introduction}

\begin{subsection}{Character stability and the AMNM property}
A fundamental question in various branches of mathematics is to determine whether ``locally approximate versions'' of a given structure are small perturbations of that structure in the global sense. Many variations of this question have been studied, often under the name ``Ulam stability''; we shall not attempt a comprehensive history here,
but for recent work in this direction see e.g.~\cite{BURGER-OZAWA-THOM,DeChiffre-Ozawa-Thom,McKenney-Vignati,Kochanek,CHOI-HORVATH-LAUSTSEN}.

This article is concerned with a form of this question regarding multiplicative linear functionals 
on commutative Banach algebras.
Given a commutative Banach algebra $A$, a bounded linear functional $\psi:A\to\Cplx$ is called \dt{multiplicative} if it satisfies $\psi(ab)=\psi(a)\psi(b)$ for all $a,b\in A$.
More generally, given $\delta>0$, we say that $\psi\in A^*$ is \dt{$\delta$-multiplicative} if the bilinear map $(a,b) \mapsto \psi(a)\psi(b)-\psi(ab)$ has norm at most $\delta$.
An obvious way to obtain examples is to take a multiplicative functional $\phi$ and put $\psi=\phi+\mu$ for some $\mu\in A^*$ of suitably small norm, thought of as a perturbation of $\phi$.
The analogue of Ulam's question now becomes: do \emph{all} ``app\-roximately multiplicative'' functionals on $A$ occur as small perturbations of multiplicative functionals?

In \cite{BEJ_AMNM1}, Johnson undertook a  systematic study of this phenomenon, and coined the acronym AMNM to describe those commutative Banach algebras for which the answer to this Ulam-type question is affirmative.
Here AMNM stands for ``\underline{a}pproximately \underline{m}ultiplicative functionals are \underline{n}ear to \underline{m}ultiplicative ones''.
Many examples are studied in \cite{BEJ_AMNM1}: it is shown there that abelian ${\rm C}^*$-algebras, $L^1$-convolution algebras of locally compact abelian groups, the Banach spaces $\ell^p$ with pointwise product, and certain algebras of holomorphic functions (including the disc algebra and $\ell^1({\mathbb Z}_+)$ and $L^1({\mathbb R}_+)$), are all AMNM. On the other hand, it is also shown in \cite{BEJ_AMNM1} that the classical Volterra algebra $L^1(0,1)$ is \emph{not} AMNM. 

In general, it seems that there is no ``one-size fits all'' method for determining if the AMNM property holds for members of some particular class of commutative Banach algebras. For instance, while many familiar uniform algebras are AMNM, see \cite{Jarosz_AMNM}, it remains unknown if $H^\infty$ of the disc is AMNM; the existence of uniform algebras which are not AMNM was open for several years, but such an example was constructed in~\cite{Sidney_AMNM}.

In the present paper, we complete a project that was initiated by work of the first author \cite{YC_wtsl-amnm} and continued in recent work of the present authors in \cite{CGP-SGF, CGP-EJC}, where the AMNM property is studied for weighted $\ell^1$-convolution algebras of certain semigroups called ``semilattices''. Our results give a complete description of those semilattices for which all such weighted convolution algebras are AMNM, or conversely, all those for which there exists some weight producing a non-AMNM algebra. The precise statements are given in the next section, once we have set up the necessary terminology.
\end{subsection}

\begin{subsection}{Weighted semilattice algebras: old and new results}
In the context of semigroup theory, a
\dt{semilattice} is a commutative semigroup in which each element is idempotent. 
Such semigroups play a particularly important role in the structure theory of semigroups, see e.g.~\cite{Howie_intro}. Moreover,
investigating
convolution algebras over semilattices fits into a well-established theme of studying how structural properties of a semigroup are reflected in properties of their associated Banach algebras. For more specific motivation, see \cite[\S1]{YC_wtsl-amnm}.

In this paper, we consider weighted convolution algebras over semilattices. A \dt{weight} on $S$ is a function $\om:S\to (0,\infty)$, and we may then define the associated weighted $\ell^1$-space $\ell^1(S,\omega)$ as the set of all functions $f:S\to\Cplx$ satisfying $\sum_{s\in S} |f(s)|\om(s)<\infty$.
Here, we only deal with weights that are \dt{submultiplicative}, i.e.~satisfying  $\om(xy)\leq \om(x)\om(y)$ for all $x$ and $y$ in $S$. 
Note that for each $x\in S$, the identity $x^2=x$ and the submultiplicative condition forces $\om(x)\geq 1$.
The condition that $\om$ be submultiplicative guarantees that  $\ell^1(S,\om)$, when equipped with the natural weighted norm, becomes a Banach algebra with respect to the convolution product
\[ (f*g)(x) \defeq \sum_{(s,t)\in S\times S \colon st=x} f(s)g(t). \]

\begin{dfn}
Given a semilattice $S$ and a submultiplicative weight $\om$, we refer to $\ell^1(S,\om)$ equipped with this convolution product as the \dt{weighted convolution algebra} of the pair $(S,\om)$.
\end{dfn}

\begin{rem}
It is well known that if $S$ is a semilattice then the unweighted convolution algebra $\ell^1(S)$ is semisimple (see e.g.~\cite[\S3]{HZ_TAMS56}).
Since $S$ is a semilattice, 
any submultiplicative weight on $S$ is bounded below by~$1$. Hence $\ell^1(S,\omega)\subseteq \ell^1(S)$, and so weighted convolution algebras on $S$ are also semisimple. Thus the algebras considered in this paper may be viewed as Banach function algebras when represented on suitable carrier spaces.
\end{rem}

To our knowledge, the first detailed study of the AMNM problem for weighted semilattice algebras was in work of the first author \cite{YC_wtsl-amnm}.
The following result was first obtained in that paper; see also Remark 4.7 of \cite{CGP-SGF} for an alternate proof. 

\begin{thmalpha}[{\cite[Example 3.13 and Theorem 3.14]{YC_wtsl-amnm}}]
\label{t:old result}
Let $S$ be a semilattice that has ``finite breadth''. Then $\ell^1(S,\omega)$ is AMNM for \emph{every} submultiplicative weight $\omega$.
\end{thmalpha}

The breadth of a semilattice
takes values in $\Nat\cup\{\infty\}$, and measures the internal complexity of the natural partial order in $S$. 
The precise definition is rather technical and will be given in Definition \ref{d:breadth} below. Semilattices of finite breadth can encompass a variety of behaviour, as shown in \cite[\S\S4-5]{ADHMS_VC2}, and so Theorem \ref{t:old result} provides a large source of semisimple commutative Banach algebras with the AMNM property.

In contrast, an explicit example is constructed in \cite[Theorem 3.4]{YC_wtsl-amnm}, of a semilattice $S$ and a weight $\omega$ for which $\ell^1(S,\omega)$ is not AMNM. This $S$ necessarily has infinite breadth, which motivated the first author to ask in
\cite[\S6]{YC_wtsl-amnm} whether the converse of Theorem~\ref{t:old result} is true. That is:
\begin{quote}
if $S$ has infinite breadth, does it always admit some submultiplicative weight $\omega$ for which $\ell^1(S,\omega)$ is non-AMNM?
\end{quote}
In \cite[Theorem 4.8]{CGP-SGF} we partially answered this question,
obtaining a positive answer whenever $S$ is a subsemilattice of $\FIN(\Omega)$,
the collection of all finite subsets of a set $\Omega$ equipped with union as the semilattice operation. 
In this article, 
by building on combinatorial techniques developed in \cite{CGP-EJC}, we are able to extend the partial answer to a full answer. This is the main new result of this article, and it provides a large supply of semisimple commutative Banach algebras which \emph{do not} have the AMNM property.
\begin{thmalpha}[Converse to Theorem \ref{t:old result}]
\label{t:new result}
Let $S$ be a semilattice that has ``infinite breadth''. Then a submultiplicative weight $\omega$ can be constructed such that $\ell^1(S,\omega)$ is not AMNM.
\end{thmalpha}

Semilattices of infinite breadth should not be seen as rare or pathological.
For example, for any infinite set $\Omega$, the power set $\cP(\Omega)$ forms a semilattice of infinite breadth (with union as the binary operation).
There are many other possibilities. 
However, it turns out that there are three particular semilattices of infinite breadth, denoted by
$\TMAX$, $\TMIN$, and $\TORT$,
which play a key role when considering general semilattices of infinite breadth.
In \cite[Theorem 1.6]{CGP-EJC}, it was shown that
if $S$ is a semilattice with infinite breadth, then there is a homomorphic image of $S$ which contains a copy of either $\TMAX$, $\TMIN$ or $\TORT$.

To prove Theorem~\ref{t:new result}, we need to prove a refined version of  \cite[Theorem 1.6]{CGP-EJC}, which we obtain using a combinatorial result with a  Ramsey-theoretic flavour
stated in Theorem~\ref{t:dichotomy}. 
The results of \cite{CGP-EJC} work with union-closed set systems $\cS\subseteq\cP(\Omega)$ and implicitly use
the notion of a \dt{spread in~$\Omega$} (see Definition~\ref{def:spread}).
In this paper we pursue a deeper study of how the set system can interact with such a spread, building up to Theorem~\ref{t:dichotomy}. The result requires too many technical definitions to be stated here, but loosely speaking it says that when $\cS$ interacts with a spread, we can control the complexity in a certain technical sense via a \dt{colouring} of the spread, unless we are in a special situation with very high complexity (``shattering'').
\end{subsection}

\begin{subsection}{Organization of the paper}
In Section~\ref{s:prelim}, we provide the preliminaries and background on weighted semilattice algebras, breadth of a semilattice, and important examples of semilattices with infinite breadth.
In Theorem~\ref{structure of infinite breadth}, we provide a self-contained statement of the main result of \cite{CGP-EJC},
describing the structure of semilattices with infinite breadth in terms of the occurrence of $\TMAX$, $\TMIN$ and $\TORT$ in them.
Then in Section~\ref{sec:weight-case 1 and 2}, we give a simple construction of non-AMNM weights for semilattices which contain $\TMAX$ or $\TMIN$ in the sense of Theorem~\ref{structure of infinite breadth}.

We devote Section~\ref{s:dichotomy} to notions of shattering and colouring, and the proof of Theorem~\ref{t:dichotomy}.
Using this we can give a refinement of our structure theorem for semilattices of infinite breadth, stated in Theorem~\ref{structure of infinite breadth-refined}. We use this structure theorem in Section~\ref{sec:weight-case 3} to produce non-AMNM weights for semilattices containing $\TORT$.
The last section of the paper contains examples showing that the construction of non-AMNM weights is somewhat delicate.
\end{subsection}
\end{section}

\begin{section}{Preliminaries}\label{s:prelim}
\begin{subsection}{Weighted semilattices and propagation}
A \dt{semilattice} is a commutative semigroup $S$ satisfying $x^2=x$ for all $x\in S$. For two elements $x,y\in S$, we say that $y$ is a \dt{multiple} of $x$ or $x$ is a \dt{factor} of $y$ or $x$ \dt{divides} $y$ and write $x\divides y$ if there exists $z\in S$ such that $y=xz$; which for the semilattice $S$ is simply equivalent to  $xy=y$.
The divisibility relation provides $S$ with a standard and canonical partial order, and to signify this aspect of the relation, we sometimes write $y\preceq x$ instead of $x\divides y$. With respect to this particular partial order, $xy$ is the \dt{meet} (or~\dt{greatest lower bound}) of $x$ and~$y$; this gives an alternative, order-theoretic definition of a semilattice. 

We repeat some terminology and notation from \cite{CGP-SGF} for the reader's convenience.
A weighted semilattice is a semilattice $S$ equipped with a submultiplicative weight. 
It will be convenient for later calculations if we switch to working with \dt{log-weights}, by which we mean functions $\lm:S \to [0,\infty)$ that satisfy $\lm(xy)\leq \lm(x)+\lm(y)$ for all $x,y\in S$. Given such a $\lm$ and $L\geq 0$, we define ``level'' $L$ as 
$W_L(S,\lm)=\{x\in S \colon \lm(x)\leq L\}$.
When there is no danger of confusion we abbreviate this to $W_L$.

\begin{dfn}[Filters in semilattices]
  \label{d:filter}
  Let $S$ be a semilattice and let $F\subseteq S$. We say that $F$ is a \dt{filter in $S$} if it is non-empty and satisfies
\[ \forall\; x,y \in S\quad (xy \in F \Longleftrightarrow x,y\in F). \]
\end{dfn}
If $E\subseteq S$, let $\filgen{E}$ denote the \dt{filter-or-empty-set generated by~$E$}, i.e.~$\filgen{\emptyset}=\emptyset$, 
and $\filgen{E}$ is the intersection of all filters $X\subseteq S$ containing $E$, if $E\neq \emptyset$.
Let $E\subseteq S$ be non-empty. Note that if $x,y\in E$ and $z\succeq xy$ (that is, $z$ is a factor of $xy$) then $z\in\filgen{E}$. Moreover, every $z\in \filgen{E}$ satisfies $z\succeq x_1\dotsb x_k$ for some $x_1,\dots, x_k\in E$.

\begin{dfn}[$\FBP_C$-stability]\label{def:FBP}
Let $C\geq 0$.
For $E\subseteq S$ we define \emph{factors of binary products} of $E$ as 
\begin{equation*}
\FBP_C(E)\defeq \set{ z\in W_C \colon \hbox{there exist $x,y \in E\cap W_C$ such that $z\succeq xy$}}.
\end{equation*}
Note that $\FBP_C(\emptyset)=\emptyset$. We also define $\FBP_C^0(E)=E\cap W_C$, and for $k\geq 1$
recursively define $\FBP_C^k(E)=\FBP_C(\FBP_C^{k-1}(E))$. 
\end{dfn}

\begin{dfn}[Propagation]
\label{d:propagation}
For $z\in\filgen{E}$, let
\begin{equation*}
V_E(z) = \inf \set{ C\geq 0 \colon \exists\  n\geq 0 \mbox{ s.t. } z\in \FBP_C^{n}(E)}.
\end{equation*}
Given $L\geq 0$, we say that $(S,\lm)$ \dt{propagates at level $L$}, or \dt{has $L$-propagation}, if
\begin{equation*}
\label{eq:sup sup finite}
\sup_{\emptyset\neq E \subseteq W_L} \sup_{z\in \filgen{E}\cap W_L}V _E(z) < \infty\;.
\end{equation*}
It is convenient to  set $V_E(z)\defeq+\infty$ whenever $z\notin\filgen{E}$.
\end{dfn}

We now connect these definitions to the original AMNM problem.
In the introduction, we defined what it means for a functional on a Banach algebra $A$ to be multiplicative  or $\delta$-multiplicative. We write $\Mult(A)$ for the set of multiplicative functionals on $A$ (note that this always includes the zero functional).

\begin{dfn}[Johnson, \cite{BEJ_AMNM1}]\label{d:AMNM}
Let $A$ be a commutative Banach algebra. We say that \dt{$A$ has the AMNM property}, or that \dt{$A$ is AMNM}, if for every $\veps>0$ we can find $\delta>0$ such that every $\delta$-multiplicative $\psi\in A^*$ satisfies $\operatorname{dist}(\psi,\Mult(A))<\veps$.
\end{dfn}

\begin{thm}[{\cite[Remark 2.7 and Theorem 3.7]{CGP-SGF}}]\label{t:equivalence}
Let $(S,\om)$ be a weighted semilattice and let $\lm=\log\om$. The following conditions are equivalent. 
\begin{newnum}
\item\label{li:(i)}
$\ell^1(S,\om)$ is AMNM.

\item\label{li:(ii)}
$(S,\lm)$ has $L$-propagation for all $L\geq 0$.
\end{newnum}
\end{thm}

\end{subsection}

\begin{subsection}{Concrete semilattices  and breadth}
To construct log-weights with the desired properties, we shall switch perspective and work with ``concrete semilattices'' that arise as union-closed set systems.

\begin{dfn}\label{d:concrete semilattices}
Let  $\Omega$ be a non-empty set, and write $\cP(\Omega)$ for its power set.
%
A  \dt{union-closed set system} or \dt{concrete semilattice} on $\Omega$ is a subset $\cS\subseteq\cP(\Omega)$ which is closed under taking finite
 unions; this is clearly a semilattice, where
 set-union serves as the binary operation.
\end{dfn}

Every semilattice can be viewed as a concrete semilattice. Indeed, given a semilattice~$S$, for each $x\in S$ let $E_x \defeq S\setminus \{ y \in S \colon x \divides y\}$. It is easily checked that $E_x\cup E_y = E_{xy}$ for all $x,y\in S$. Therefore, the function $E_\blob: S \to \cP(S)$, $x\mapsto E_x$, defines an injective semilattice homomorphism from $S$ into $(\cP(S),\cup)$. This is sometimes known as the \dt{Cayley embedding} of a semilattice. Using this result, from this point onward, we assume every semilattice is a union-closed set system. (However, see Remark \ref{r:different-reps} for some subtleties.)

We shall use the following notational conventions when working with set systems on a given set $\Omega$.
Elements of $\Omega$ will usually be denoted by lower-case Greek letters.
Set systems on $\Omega$ (i.e.~subsets of $\cP(\Omega)$) will usually be denoted by letters such as $\cB$, $\cS$, etc. If $\cB$ is such a set system, we refer to \dt{members of $\cB$} rather than elements. If $\cB$ and $\cS$ are set systems on $\Omega$ we denote their union and intersection by $\cB\sscup\cS$ and $\cB\sscap \cS$.
Members of a set system $\cS$ will be denoted by letters such as $\ssa$, $\ssb$, $\ssp$, etc., and we write $\ssa\cup\ssb$ and $\ssa\cap \ssb$ for their union and intersection respectively.
If it happens that $\ssa$ and $\ssb$ are disjoint subsets of $\Omega$ we shall sometimes emphasise this by writing their union as $\ssa\dotcup\ssb$.
Note that in a concrete semilattice $\cS$, $\ssa$ is a factor of $\ssb$ precisely when $\ssa\subseteq \ssb$.

%

The following definitions are notions from lattice theory, restated in the current setting of union-closed set systems.

\begin{dfn}\label{d:join}
Given $\cF\subseteq\cP(\Omega)$, the \dt{join of $\cF$} is the set $\join(\cF)\defeq \bigsetcup_{\ssx\in\cF} \ssx$. If $\cF$ is a finite subset of a union-closed set system $\cS$, then $\join(\cF)\in\cS$.
\end{dfn}

\begin{dfn}\label{d:breadth}
Let $\cS$ be a union-closed set system. Given a finite, non-empty subset $\cE\subseteq \cS$, we say $\cE$ is \dt{compressible} if there exists a proper subset $\cE'\subset \cE$ 
such that $\join(\cE')=\join(\cE)$; otherwise, we say $\cE$ is \dt{incompressible}.
The \dt{breadth} of a semilattice $\cS$ is defined to be
\begin{align*} 
b(\cS) &=\inf\{ n\in\Nat \colon \text{every $E\subseteq S$ with $n+1$ members is compressible} \} \\
&=\sup\{n\in\Nat \colon \text{$\cS$ has an incompressible subset with $n$ elements} \}. 
\end{align*}
\end{dfn}

\begin{eg}
We write $\FIN(\Omega)$ for the set of all finite subsets of $\Omega$; this is a concrete semilattice on $\Omega$. Note that when $\Omega$ is infinite, $\FIN(\Omega)$ has infinite breadth, since for any $\gamma_1,\dots, \gamma_n\in \Omega$ the set $\{ \{\gamma_j\} \colon 1\leq j\leq n\}$ is an incompressible subset of $\FIN(\Omega)$. 
\end{eg}

A little thought shows that the notion of compressibility for a subset of $\cS$ only depends on the underlying semigroup structure of~$(\cS,\cup)$. Therefore, the breadth of a semilattice $S$ is an \emph{intrinsic} invariant, which does not depend on any particular concrete representation of $S$ as a union-closed set system. (For a direct definition without using concrete representations, see e.g.\ \cite[Definition 4.5]{CGP-SGF}.)

The breadth of a semilattice sheds some light on its structure, and is related to more familiar order-theoretic concepts such as height and width. (Some basic links, with references, are surveyed in \cite[Section 4.1]{ADHMS_VC2}.) For instance, by examining incompressible subsets, one sees that if $b(\cS)\geq n$ then $\cS$ contains a chain (totally ordered subset) and an antichain (subset in which no two elements are comparable) both of cardinality~$n$. In particular, a semilattice $\cS$ has breadth~$1$ exactly when the poset $(\cS,\preceq)$ is totally ordered.

\begin{subsection}{Special set systems}
If $b(\cS)=\infty$, then there are arbitrarily large finite subsets of $\cS$ that are incompressible. However, unlike the example of $\FIN(\Omega)$, we cannot always arrange for these to be nested in an infinite sequence $\cF_1\subseteq \cF_2\subseteq \dots$
This can be seen very clearly with the three key examples that will be introduced in Definition \ref{d:three key examples}.
\end{subsection}

Before defining these examples, we introduce some terminology that will be convenient for Section \ref{s:dichotomy}.

\begin{dfn}\label{def:spread}
A \dt{spread} in  a set $\Omega$ is a sequence $\cE=(E_n)_{n\geq 1}$ of finite non-empty subsets of $\Omega$ which are pairwise disjoint and satisfy $|E_n | \to \infty$.  
 A \dt{refinement of $\cE$} is a spread $\cF=(F_j)_{j\geq 1}$ with the property that each $F_j$ is contained in some $E_{n(j)}$ and $n(j)\neq n(k)$ whenever $j\neq k$.
\end{dfn}
By a slight abuse of notation: if $\cE=(E_n)_{n\geq 1}$ is a spread, we shall write $\join(\cE)$ for $\bigcup_{n\geq 1} E_n$.

\begin{dfn}[Three special set systems]
\label{d:three key examples}
Let $\cE=(E_n)_{n\geq 1}$ be a spread in  a set $\Omega$. For $n\in\Nat$, let $E_{<n}\defeq E_1\dotcup \dots \dotcup E_{n-1}$ (with the convention that $E_{<1} =\emptyset$) and let $E_{>n} \defeq \dot{\bigsetcup}_{j\geq n+1} E_j$. We now define the following set systems on $\Omega$:
\[
\begin{aligned}
\TMAX(\cE) & \defeq \bigsscup_{n\geq 1} \bigsscup_{\emptyset\neq \ssa\subseteq E_n} \left\{ E_{<n}\dotcup \ssa \right\}, 
\\
\TMIN(\cE) & \defeq \bigsscup_{n\geq 1} \bigsscup_{\emptyset\neq \ssa\subseteq E_n} \set{ \ssa\dotcup E_{>n}},  
\\
\TORT(\cE) & \defeq \bigsscup_{n\geq 1} \bigsscup_{\emptyset\neq \ssa\subseteq E_n} \left\{ E_{<n}\dotcup \ssa\dotcup E_{>n} \right\}.
\end{aligned}
\]
\end{dfn}

Note that this definition is slightly more general than \cite[Definition 1.3]{CGP-EJC}, where $|E_{n}|=n+1$ was assumed. 

\begin{figure}[hpt]
\centering
\includegraphics[scale=0.8]{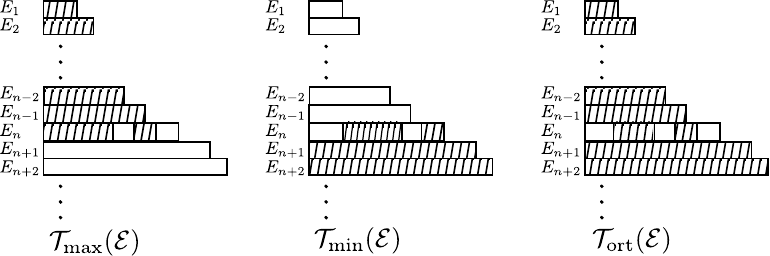}
\caption{Typical members of $\TMAX$, $\TMIN$, and $\TORT$ at level $n$}
\label{fig1}
\end{figure}

The key role of these examples was shown in \cite{CGP-EJC} by the following result.

\begin{thm}[{\cite[Theorem 1.6]{CGP-EJC}}]\label{structure of infinite breadth}
Let $\cS$ be a union-closed set system on $\Omega$. Suppose  $\cS$ has  infinite breadth.
Then there is a spread in $\Omega$, denoted by $\cE=(E_n)_{n\in{\mathbb N}}$, such that at least one of the following statements holds:
\begin{newnum}
\item\label{case:TMAX}
 $\left\{\ssx\cap \join(\cE) \colon \ssx\in \cS\right\}\supseteq \TMAX(\cE)$.
\item\label{case:TMIN}
 $\left\{\ssx\cap \join(\cE) \colon \ssx\in \cS\right\}\supseteq\TMIN(\cE)$.
\item\label{case:TORT}
 $\left\{\ssx\cap \join(\cE) \colon \ssx\in \cS\right\}\supseteq \TORT(\cE)$.
\end{newnum}
\end{thm}

\begin{rem}\label{r:different-reps}
In applying this theorem to prove Theorem~\ref{t:new result}, we are choosing to represent an abstract semilattice $S$ as a concrete semilattice $\cS$. The following example shows that the choice of ``concrete representation'' may affect which case of Theorem \ref{structure of infinite breadth} occurs. (It plays no role in the proof of Theorem \ref{t:new result}, and can be skipped on a first reading.)
\end{rem}

\begin{eg}\label{nonpersistence of TMIN and TORT}
Set $E_n\defeq\set{(n,k) \colon 1\le k\le n} \subset \Nat^2$, and set
\[
	\cE\defeq(E_n)_{n\ge 1},\quad \Omega_0\defeq\join(\cE),\quad\text{and}\quad \Omega\defeq\Omega_0\sqcup\Nat.
\]
Here $\sqcup$ denotes the formal disjoint union of sets.
For each $\ssa\in\TMIN(\cE)$, define $\level(\ssa)$ to be the level of $\ssa$ as indicated in Figure \ref{fig1}.  
Define 
\begin{equation}\label{eq:tmin not tort}
	\cS\defeq\set{\ssa\sqcup\set{1,\ldots,m} \colon\ \ssa\in\TMIN(\cE),\ m\in\Nat,\ m\ge \level(\ssa)\ge 2}.
\end{equation}

For $\ssa_i\in\TMIN(\cE)$ and $m_i \geq \level(\ssa_i)$ ($i=1,2$),
\[
(\ssa_1 \sqcup \set{1,\ldots,m_1}) \cup
(\ssa_2 \sqcup \set{1,\ldots,m_2})
= (\ssa_1\cup\ssa_2) \sqcup \set{1,\ldots, \max(m_1,m_2)}
\]
where $\level(\ssa_1\cup\ssa_2)=\min( \level(\ssa_1),\level(\ssa_2))\leq \max(m_1,m_2)$. 
Thus $\cS$ is a concrete semilattice on $\Omega$. It satisfies \ref{case:TMIN} of Theorem \ref{structure of infinite breadth}, since $\{\ssx \cap \join(\cE) \colon \ssx\in\cS\} = \TMIN(\cE)$ by construction.

The following short argument shows that $\cS$ satisfies neither \ref{case:TMAX} or \ref{case:TORT} of that theorem. Suppose that $\cF=(F_n)_{n=1}^\infty$ is a spread on $\Omega$ such that $\left\{\ssx\cap \join(\cF) \colon \ssx\in \cS\right\}$ contains either $\TMAX(\cF)$ or~$\TORT(\cF)$.
We may suppose that $\join(\cF)\subseteq \Omega_0$ (if not, simply replace each $F_n$ by $F_n\cap\Omega_0$ and remove any resulting empty set, noting that $\abs{F_n\setminus\Omega_0}\le 1$ for all~$n$).
Then there exist only finitely many members of $\left\{\ssx\cap \join(\cF) \colon \ssx\in \cS\right\}$ that meet $F_1$, by~\eqref{eq:tmin not tort}. However, there are infinitely many members of  $\TMAX(\cF)$ and of $\TORT(\cF)$ that meet $F_1$; a contradiction.

On the other hand, set $\Omega'\defeq\Omega_0\times \Nat$, and define
\begin{equation}\label{eq:tort not tmin}
\begin{aligned}
	\cS'\defeq\set{(\ssa\times\Nat)\cup\left(\Omega_0\times\set{1,\ldots,m}\right) \colon\ \ssa\in\TMIN(\cE),\ m\in\Nat,\ m\ge \level(\ssa)\ge 2}.
\end{aligned}
\end{equation}
Then $\cS'$ is a concrete semilattice on $\Omega'$, and a little thought shows that $\cS$ and $\cS'$ are isomorphic as abstract semilattices. However, we claim that it does not satisfy \ref{case:TMIN}.
\end{eg}

\begin{proof}[Proof of the claim]
First of all, observe that for any sequence $(\ssx_j)_{j\ge 1}$ of distinct members of $\cS'$ one has $\bigcup_{j=1}^\infty \ssx_j=\Omega'$. 
Indeed, write $\ssx_j=(\ssa_j\times\Nat)\cup\left(\Omega_0\times\set{1,\ldots,m_j}\right)$. Since $\{\ssx_j\}$ is infinite, \eqref{eq:tort not tmin} implies that $\{m_j\}$ is not bounded, and so $\bigcup_{j=1}^\infty \left(\Omega_0\times\set{1,\ldots,m_j}\right)=\Omega'$.

Now suppose there exists a spread $\cF=(F_k)_{k\ge 1}$ such that
$\{ \ssx\cap \join(\cF) \colon \ssx\in \cS'\} \supseteq\TMIN(\cF)$.
For each $j\in\Nat$, choose $\ssx_j\in \cS'$ such that $\ssx_j\cap \join(\cF)$ belongs to level $j+2$ in $\TMIN(\cF)$. Then 
\begin{align*}
        \join(\cF)\neq \bigcup_{j=1}^\infty\left(\ssx_j\cap \join(\cF)\right)=\left(\bigcup_{j=1}^\infty\ssx_j\right)\cap \join(\cF)=\Omega'\cap\join(\cF)=\join(\cF)
\end{align*}
a contradiction.
\end{proof}

By similar reasoning, one can give a direct proof that $(\cS',\Omega')$ does not satisfy \ref{case:TMAX}. We omit the details, since it also follows from the next proposition, which shows that the $\TMAX(\cE)$ case of Theorem \ref{structure of infinite breadth} is better behaved with respect to isomorphism of abstract semilattices.
Finally, by Theorem~\ref{structure of infinite breadth} the pair  $(\cS',\Omega')$  satisfies \ref{case:TORT} of Theorem \ref{structure of infinite breadth}.

\begin{prop}\label{persistence of TMAX}
Let $S$ be a semilattice and let $\Omega$, $\Omega'$ be infinite sets with injective homomorphisms
$\Theta: S\to (\cP(\Omega),\cup)$ and $\Theta': \cS\to (\cP(\Omega'),\cup)$. Suppose there is  a spread $\cE$ in $\Omega$ such that 
$\{\Theta(s)\cap \join(\cE): s\in S\}\supseteq \TMAX(\cE)$. Then there is a spread $\cE'$ in $\Omega'$ such that 
$\{\Theta'(s)\cap \join(\cE'): s\in S\}\supseteq \TMAX(\cE')$.
\end{prop}

\begin{proof}
For simplicity of notation,
 we shall write $a$ or $b_n$ for elements of $S$, while writing $\ssa$ or $\ssb_n$ for the corresponding members of~$\cS$ and  $\ssa'$ or $\ssb'_n$ for the corresponding members of~$\cS'$. Also, without loss of generality, we suppose that $\cE=(E_n)_{n\geq 1}$ with $|E_n |=n$ for all $n$; say $E_n=\set{\gamma_{nj} \colon 1\le j\le n}$.
 
 For $1\le j\le n$, let $a_{nj}$ be an element of $S$ such that $\ssa_{nj}\cap \join(\cE)$ is a member of $\TMAX(\cE)$ that meets $E_n$ at the singleton $\set{\gamma_{nj}}$. Then, for each $n\in\Nat$, $a_{nj}$ does not divide $a_{ml}a_{nk}$ where $m<n$ and $k\neq j$, and so we can find an element $\gamma'_{nj}$ of $\Omega'$ that belongs to $\ssa'_{nj}$ but not to any of $\ssa'_{ml}$ and $\ssa'_{nk}$ where $m<n$ and $k\neq j$. Define $E'_n\defeq\set{\gamma'_{nj} \colon 1\le j\le n}$, and then $\cE'\defeq(E'_n)_{n\geq 1}$. Then $\cE'$ is a spread in $\Omega'$, and it satisfies
\[
 \{\ssx\cap\join(\cE') \colon \ssx\in\cS'\} \supseteq\TMAX(\cE'),
\]
since for each $1\le j\le n$ we have
\[
	\left(\bigcup_{1\le l\le m<n}\ssa'_{ml}\cup \ssa'_{nj}\right)\cap \join(\cE')=E'_{<n}\cup\set{\gamma'_{nj}}.
\]
\end{proof}

\end{subsection}

\end{section}
\begin{section}{Constructing non-AMNM weights in the $\TMAX$ and $\TMIN$ cases}\label{sec:weight-case 1 and 2}
\begin{prop}\label{prop:tmax-case}
Let $\cS$ be a union-closed set system on $\Omega$ and suppose $\cE$ is a spread in $\Omega$ such that Case~\ref{case:TMAX} of Theorem~\ref{structure of infinite breadth} holds. 
Then there is a log-weight on $\cS$ which fails to propagate at the first level.
\end{prop}
\begin{proof}
Suppose that there is a spread $\cE=(E_n)_{n\geq 1}$ in $\Omega$ such that $\left\{\ssx\cap \join(\cE) \colon \ssx\in \cS\right\}$ contains $\TMAX(\cE)$.
Passing to a refinement of $\cE$ if necessary, we may suppose that $\abs{E_n}=n+1$ for all~$n$ (this just simplifies some of the following formulas and arguments).

 For each $\ssx\in \cP(\Omega)$, define
\begin{align}\label{variant eq:MG-weight type I}
	\lm(\ssx)\defeq\begin{cases}
		0&\ \textrm{if there are \emph{no} or \emph{infinitely many} $n$ such that $E_n\cap \ssx\neq\emptyset$,}\\
		0&\ \textrm{if $ E_n\subseteq \ssx$, and $E_k\cap \ssx=\emptyset$ for $k>n$,}\\ 
		\abs{\ssx\cap E_n}&\ \textrm{if $E_n\cap\ssx\neq \emptyset$, $E_n \not\subseteq \ssx$, and $E_k\cap \ssx=\emptyset$ for $k>n$}.
	\end{cases}
\end{align}

\begin{claim}\label{variant p:MG-subadd type I}
Let $\ssx,\ssy\in \cP(\Omega)$. Then $\lm(\ssx\setcup\ssy)\leq \lm(\ssx)+\lm(\ssy)$.
\end{claim}

\begin{proof}[Proof of the claim.]
Fix $\ssx,\ssy\in \cP(\Omega)$. The desired inequality is trivial
when $\ssx\cap \join(\cE)=\emptyset$ or when $\ssy\cap \join(\cE)=\emptyset$.
On the other hand: if there are infinitely many $n$ such that $E_n\cap \ssx\neq\emptyset$, or infinitely many $n$ such that $E_n\cap \ssy\neq\emptyset$, then the same is true for $\ssx\setcup\ssy$, and so the inequality follows, since both $\lm(\ssy)\ge 0$ and $\lm(\ssx)\ge 0$.

Otherwise, let $m$ be the largest natural number such that $ E_m\cap \ssx\neq\emptyset$, and let $n$ be the corresponding number for $\ssy$. If $m\neq n$, then without loss of generality we suppose $m> n$. Then $m$ is also the largest number such that $ E_m\cap (\ssx\setcup\ssy)\neq\emptyset$. Moreover $E_m\cap (\ssx\setcup\ssy)= E_m\cap\ssx$, so the inequality follows, since $\lm(\ssy)\ge 0$.

If not, $m=n$ is the largest number such that $E_m\cap (\ssx\setcup\ssy)\neq\emptyset$. If furthermore either $E_m\subseteq \ssx$ or $E_m\subseteq \ssy$, then $E_m\subseteq\ssx\setcup\ssy$, and the inequality is again obvious. Otherwise, we see that
\[
	\lm(\ssx\setcup \ssy)\le\abs{(\ssx\setcup \ssy)\cap E_m}\le \abs{\ssx\cap E_m}+\abs{\ssy\cap E_m}=\lm(\ssx)+\lm(\ssy)\,.
\]
This completes the proof of our claim.
\end{proof}

We will now show that $(\cS,\lm)$ does not have $1$-propagation.
Let $n\in\Nat$.
Since $\TMAX(\cE)$ is contained in $\left\{\ssx\cap \join(\cE) \colon \ssx\in \cS\right\}$, there is a size-$(n+1)$ subset $\cF_n\subseteq \cS$ such that:
\begin{itemize}
\item  for each $\ssa\in\cF_n$,  $E_j\subseteq\ssa$ for $j<n$, $\ssa\cap E_n$ is a singleton, and $\ssa\cap E_j=\emptyset$ for $j>n$;
\item $\ssa\cap E_n$ gives different singletons for different $\ssa\in\cF_n$.
\end{itemize}
From our construction, $\lm(\ssa)=1$ for every $\ssa\in\cF_n$.

Set $\ssb_n\defeq\join(\cF_n)$. Then $\lm(\ssb_n)=0$, and so $\ssb_n\in \filgen{\cF_n}\sscap W_1$. To finish the proof it suffices to show that $V_{\cF_n}(\ssb_n)\to \infty$ as $n\to\infty$, which we do as follows. Let $C\geq 1$ be such that $\ssb_n\in\bigcup_{k=0}^\infty \FBP_C^k(\cF_n)$.  Let $m\geq 1$ be minimal  with respect to the following property:
\begin{quote}
there exists $\ssa\in \FBP_C^m(\cF_n)$ such that  $E_n\subseteq \ssa$.
\end{quote}
(Such an $m$ exists by our assumption, since $ E_n\subseteq \ssb_n$.) By minimality there are $\ssa_1$ and $\ssa_2$ in $\FBP_C^{m-1}(\cF_n)$ such that $ E_n$ is contained in $\ssa_1\setcup\ssa_2$, yet $ E_n\not\subseteq\ssa_1$ and $ E_n\not\subseteq\ssa_2$. (When $m=1$, our convention here is that $\FBP_C^0(\cF_n)=\cF_n$ and, since $\abs{E_n}\ge 2$, $E_n$ cannot be a subset of any member of $\cF_n$ either.)

Let $i\in\{1,2\}$. By the previous remarks, $\ssa_i\cap E_n$ is a proper, nonempty subset of~$E_n$. By Definition~\ref{def:FBP}, $\ssa_i \subseteq\ssb_n$. Thus $n$ is the largest natural number $k$ such that $\ssa_i\cap E_k\neq\emptyset$, and so $\lm(\ssa_i)=\abs{\ssa_i\cap E_n}$. Hence $\abs{ E_n} \leq \lm(\ssa_1)+\lm(\ssa_2)\leq 2C$, with the last inequality following because $\ssa_1,\ssa_2\in W_C$. Therefore, $V_{\cF_n}(\ssb_n) \geq \frac{1}{2}\abs{E_n} \to \infty$, completing the proof.  
\end{proof}

\begin{prop}\label{prop:tmin-case}
Let $\cS$ be a union-closed set system on $\Omega$ and suppose $\cE$ is a spread in $\Omega$ such that Case~\ref{case:TMIN} of Theorem~\ref{structure of infinite breadth} holds. 
Then there is a log-weight on $\cS$ which fails to propagate at the first level.
\end{prop}
\begin{proof}
Suppose that there is a spread $\cE=(E_n)_{n\geq 1}$ in $\Omega$ such that Case~\ref{case:TMIN} of Theorem~\ref{structure of infinite breadth} holds. Then, for each $\ssx\in \cP(\Omega)$, define
\begin{align}\label{variant eq:MG-weight type II}
\lm(\ssx)\defeq\begin{cases}
0&\ \textrm{$\ssx\cap \join(\cE)=\emptyset$},\\
0&\ \textrm{if $ E_n\subseteq \ssx$, and $E_k\cap \ssx=\emptyset$ for $1\leq k< n$,}\\ 
\abs{\ssx\cap E_{n}}&\ \textrm{if $E_n\cap\ssx\neq \emptyset$, $E_n\not\subseteq \ssx$, and $E_k\cap \ssx=\emptyset$ for $1\leq k< n$}.
\end{cases}
\end{align}

\begin{claim}
In this case, $\lm$ is a log-weight on $\cS$, and $(\cS,\lm)$ does not have $1$-propagation.
\end{claim}
The rest of the proof is similar to that of Proposition~\ref{prop:tmax-case}, but slightly easier, so we omit the details.
\end{proof}
\end{section}

\begin{section}{Shattering, colouring, and a Ramsey-theoretic result}
\label{s:dichotomy}
Constructing non-AMNM weights, when Case~\ref{case:TORT} of Theorem~\ref{structure of infinite breadth} holds, is much more involved. In this case, we need to gain a deeper understanding of the incompressible subsets of the semilattice.

\begin{dfn}[Shattering a spread]
Let $\cE=(E_n)_{n\geq 1}$ be a spread in $\Omega$ and let $(\ssa_j)_{j\geq 1}$ be a sequence of subsets of $\Omega$. We say that this sequence \dt{shatters $\cE$} if, for every $m\in\Nat$ and every $m$-tuple $(\ssy_1,\dots,\ssy_m)$ such that $\ssy_j\in\{\ssa_j,{\ssa_j}^c\}$ for $j=1,\dots, m$,
\[ \lim_{n\to\infty} \left\vert E_n \setcap \bigsetcap\nolimits_{j=1}^m \ssy_j \right\vert = \infty.\]
\end{dfn}

It is easy to verify that the if a sequence $(\ssa_j)_{j\geq 1}$ shatters $\cE$, then for every $\ssy_j\in \{\ssa_j,{\ssa_j}^c\}$ and every $1\leq j\leq m$, we must have 
$\ssy_1\cap \dots\cap \ssy_m \neq \emptyset$.
As a consequence, the sets $\ssa_1,\dots,\ssa_m$ are mutually distinct and form an incompressible subset of $\cP(\Omega)$. So if we start with a shattering sequence $(\ssa_j)$ and allow ourselves to take complements and binary unions (i.e.~we generate a ring of sets) then the resulting collection has high complexity; moreover, this complexity is seen inside each $E_n$ once we take $n$ sufficiently large. 

\begin{dfn}
Given a spread $\cE=(E_n)_{n\geq 1}$, a partition of $\Omega$ into finitely many disjoint subsets $\Omega= C_0 \dotcup \dots \dotcup C_d$ is said to \dt{colour $\cE$} if $\lim_n |C_j\setcap E_n | =\infty$ for each $j=0,\ldots, d$. We call $\cC=\{C_0, \dots, C_d\}$ a \dt{colouring of $\cE$.} 
\end{dfn}
\begin{dfn}[Decisive colourings]
Let $\cS\subseteq \cP(\Omega)$ be a set system, let $\cE$ be a spread in $\Omega$ and let $\cC$ be a coloring of $\cE$.  We say that $\cC$ \dt{decides $\cS$ with respect to $\cE$}, or is \dt{$\cS$-decisive with respect to $\cE$}, if there exists a colour class $C_0\in\cC$ such that for every $\ssx\in\cS$ satisfies
\begin{equation}\label{eq:decisive}
\sup_{n\geq 1}\ \min\{ | \ssx \setcap C_0 \setcap E_n | \;;\; | \ssx^c \setcap C \setcap E_n |, C\in \cC \} <\infty .
\end{equation}
Such a  $C_0$ is said to be a \dt{decisive colour class} (for $\cS$ with respect to $\cE$).  
\end{dfn}

Informally speaking, when we have a decisive colour class $C_0$, each $\ssx\in\cS$ must either have small intersection with $C_0\setcap E_n$, or else have large intersection with $C\setcap E_n$ for some $C\in\cC$,  once $n$ is sufficiently large. 

The following theorem tells us, loosely speaking, that unless we are in the highly fragmented case, we can find a decisive colouring.

\begin{thm}\label{t:dichotomy}
Let $\cS\subseteq \cP(\Omega)$ be a union-closed set system, and let $\cE$ be a spread in $\Omega$. Then at least one of the following conclusions holds:
\begin{itemize}
\item[{\rm(S)}] there is a spread $\cG$ that refines $\cE$, and a sequence $(\ssa_i)_{i\geq 1}$ in $\cS$ which shatters $\cG$;
\item[{\rm(D)}] there is a spread $\cF$ that refines $\cE$, and an $\cS$-decisive colouring $\cC$ of $\cF$.
\end{itemize}
\end{thm}
The proof of Theorem \ref{t:dichotomy}
 requires a double induction, and we isolate part of it as a separate lemma. The following terminology is introduced to streamline the presentation. Let $\cE=(E_n)_{n\geq 1}$ be a spread, and let $D$ and $F$ be subsets of $\Omega$. We say that $F$ \dt{halves $D$ with respect to $\cE$}, if
\[
	\lim_{n\to\infty} \abs{D\cap F\cap E_n}=\lim_{n\to\infty} \abs{(D\setminus F)\cap E_n}=\infty.
\]
Also, if $N\subseteq \Nat$ is an infinite subset, and $(t_n)_{n\geq 1}$ is a sequence in $[0,\infty)$, we say that \dt{$t_n\to \infty$ along $n\in N$} if $\lim_{k\to\infty} t_{n_k} = \infty$, where $N=\{n_k \colon k\in\Nat\}$.

\begin{lem}\label{l:inductive-step}
Let $\cS\subseteq\cP(\Omega)$ be a  union-closed set system. Let $\cE=(E_n)_{n\geq 1}$ be a spread in $\Omega$, and let $\cC$ be a colouring of $\cE$. Then at least one of the following conclusions holds:
\begin{newnum}
\item\label{li:Case1} there is an infinite set $N\subseteq \Nat$ and some $\ssy\in \cS$  which halves every $C\in \cC$ with respect to $(E_n)_{n\in N}$;
\item\label{li:Case2} there is a spread $\cF$ refining $\cE$, such that $\cC$ is an $\cS$-decisive colouring of $\cF$.
\end{newnum}
\end{lem}

\begin{proof}
The idea is as follows: we attempt to construct, by iteration, members of $\cS$ that are closer and closer to satisfying the property in Case~\ref{li:Case1}. At each stage of the iteration we will be able to continue, unless we find ourselves in Case~\ref{li:Case2}. Therefore, if Case~\ref{li:Case2} does not hold, our iteration will run successfully until Case~\ref{li:Case1} is satisfied.

Assume from now on that Case~\ref{li:Case2} does not hold, and enumerate the members of $\cC$ as $C_1,\dots, C_d$. Since $\cC$ is not an $\cS$-decisive colouring of $\cE$,  there exists $\ssy_1\in\cS$ such that
\[ \sup_n \min\{ |\ssy_1 \setcap C_1\setcap E_n |\;;\; |{\ssy_1}^c\setcap C_j\setcap E_n | \;,\;1\leq j\leq d\} = \infty. \]
Passing to an appropriate subsequence, there exists an infinite $N_1\subseteq \Nat$ such that:
\begin{itemize}
\item $\abs{\ssy_1\setcap C_1\setcap E_n}\to \infty$ along $n\in N_1$; and
\item $| {\ssy_1}^c\setcap C_j\setcap E_n| \to\infty$ along $n\in N_1$ for every $1\leq j\leq d$.
\end{itemize}
Now let $2\leq k \leq d$. Suppose there are $\ssy_{k-1}\in \cS$ and an infinite $N_{k-1}\subseteq\Nat$ such that:
\begin{equation}\label{eq:Cond1}
| \ssy_{k-1} \setcap C_i \setcap E_n | \to\infty \text{ along $n\in N_{k-1}$ for all $1\leq i\leq k-1$}
\end{equation}
and
\begin{equation}\label{eq:Cond2}
| {\ssy_{k-1}}^c \setcap C_j \setcap E_n | \to \infty \text{ along $n\in N_{k-1}$ for all $1\leq j \leq d$.}
\end{equation}
Informally, \eqref{eq:Cond1} says that $\ssy_{k-1}\setcap C_i$ is ``not too sparse'' relative to the spread $(E_n)_{n\in N_{k-1}}$, for $1\leq i\leq k-1$, while \eqref{eq:Cond2} says that $\ssy_{k-1}\setcap C_j$ is ``not too dense'' relative to the same spread, for all~$j$.

We might have ${\ssy_{k-1}}^c\cap E_m=\emptyset$ for some $m\in N_{k-1}$. However, by \eqref{eq:Cond2} we can assume (after replacing $N_{k-1}$ with some cofinal subset, if necessary) that the sequence $({\ssy_{k-1}}^c\setcap E_n)_{n\in N_{k-1}}$ is a spread in $\Omega$, which we denote by $\cE_k$. By construction, $\cE_k$ refines $\cE$, and by condition \eqref{eq:Cond2} again, $\cC$  colours $\cE_k$. Since this is not an $\cS$-decisive colouring, in particular $C_k$ must be indecisive. By the same reasoning as before, there is some $\ssx\in \cS$ and an infinite $N_k\subseteq N_{k-1}$ such that:
\begin{equation}\label{eq:Cond3}
|\ssx \setcap C_k \setcap ({\ssy_{k-1}}^c \setcap E_n) | \to \infty \text{ along $n\in N_k$}
\end{equation}
and
\begin{equation}\label{eq:Cond4}
|\ssx^c \setcap C_j \setcap ({\ssy_{k-1}}^c \setcap E_n) | \to \infty \text{ along $n\in N_k$, for every $1\leq j \leq d$.}
\end{equation}

Let $\ssy_k\defeq \ssy_{k-1}\setcup \ssx$, which belongs to $\cS$ since $\cS$ is union-closed. Since $N_k\subseteq N_{k-1}$ and $\ssy_k\supseteq \ssy_{k-1}$, condition~\eqref{eq:Cond1} implies
\[ | \ssy_k\setcap C_i \setcap E_n | \to \infty \text{ along $n\in N_k$, for all $1\leq i\leq k-1$} \;;\]
and since $\ssy_k \supseteq \ssx$, condition~\eqref{eq:Cond3} implies
\[ | \ssy_k\setcap C_k \setcap E_n | \to \infty \text{ along $n\in N_k$} \;.\]
But by the definition of $\ssy_k$, condition \eqref{eq:Cond4} can be rephrased as
\[ |{\ssy_k}^c \setcap C_j \setcap E_n |  \to \infty \text{ along $n\in N_k$, for every $1\leq j\leq d$} \;.\]
Thus the induction can continue. We end up with $\ssy_d\in\cS$ and an infinite subset $N_d\subseteq \Nat$ such that $\ssy_d$ halves $C_j$ with respect to $(E_n)_{n\in N_d}$ for every $1\leq j\leq d$, i.e.~we are in Case~\ref{li:Case1} of the lemma.
\end{proof}

\begin{proof}[Proof of Theorem~\ref{t:dichotomy}]
Suppose that Case (D) of the theorem does not hold. The following notation will be useful: given $\ssx_1,\dots,\ssx_m\in \cP(\Omega)$, consider
\[ \Gamma(\ssx_1,\dots,\ssx_m) \defeq \left\{ \bigsetcap_{j=1}^m {\ssy_j} \; \colon\; \ssy_j\in\{\ssx_j,\ssx_j^c\}\;\text{for each $j=1,\ldots,m$}\right\}. \]
This is a partition of~$\Omega$, although some members of the partition might be empty.

Now let $\cE_0=\cE$ and let $\cC_0$ denote the trivial colouring $\{\Omega\}$.  Apply Lemma~\ref{l:inductive-step} to the pair $(\cE_0,\cC_0)$. Case~\ref{li:Case2} of the lemma does not hold, since otherwise we would be in Case (D) of the theorem. Hence we are in Case~\ref{li:Case1} of the lemma, so there exists an infinite set $N_1\subseteq\Nat$ and some $\ssa_1\in \cS$ which halves $\Omega$ with respect to $(E_n)_{n\in N_1}$.

Suppose that for some $k\geq 1$, we have found $\ssa_1,\dots,\ssa_k\in\cS$ and an infinite subset $N_k\subseteq \Nat$, such that $\Gamma(\ssa_1,\dots,\ssa_k)$ colours the spread $(E_n)_{n\in N_k}$. This colouring cannot be $\cS$-decisive  (otherwise we would be in Case (D) of the theorem, contrary to assumption). Hence, by Lemma~\ref{l:inductive-step} there exist some infinite $N_{k+1} \subseteq N_k$ and some $\ssa_{k+1}\in\cS$, such that $\ssa_{k+1}$ halves $C$ with respect to $(E_n)_{n\in N_{k+1}}$ for each $C\in\Gamma(\ssa_1,\dots,\ssa_k)$. Now $\Gamma(\ssa_1,\dots,\ssa_{k+1})$ is a colouring of the spread $(E_n)_{n\in N_{k+1}}$.

Continuing in this way, we inductively construct a sequence $(\ssa_n)_{n\geq 1}$ in $\cS$, and a descending chain of infinite subsets of $\Nat$, $N_1\supseteq N_2 \supseteq \dots$, such that:
\begin{quote}
for each $m\geq 1$ and each $C\in \Gamma(\ssa_1,\dots,\ssa_m)$, $| C\setcap E_n | \to \infty$ along $n\in N_m$.
\end{quote}

Since $N_1\supseteq N_2\supseteq \ldots$ is a decreasing sequence of infinite subsets of $\Nat$, we can extract a diagonal subsequence $n(1)<n(2)<n(3)<\dots$ satisfying $n(k)\in N_j$ for every $j\leq k$. For each~$m$, $(n(i))_{i\geq m}$ is a subsequence of $N_m$, and so for each $C\in \Gamma(\ssa_1,\dots,\ssa_m)$ we have
\[ \lim_{i\to\infty} |C\setcap E_{n(i)}| =\lim_{N_m \ni n\to\infty} |C\setcap E_n| =\infty\;.\]
 Therefore, if we define $G_i = E_{n(i)}$, the sequence $(G_i)_{i\geq 1}$ is a spread which is shattered by the sequence $(\ssa_j)_{j\geq 1}$, and we are in Case~(S) as required.
\end{proof}

We can now refine the structure result in Theorem~\ref{structure of infinite breadth}. This is the bare minimum that we need  for our weight construction in the next section.
\begin{thm}[Refined version of Theorem~\ref{structure of infinite breadth}]
\label{structure of infinite breadth-refined}
Let  $\cS$ be a union-closed set system on $\Omega$. If $\cS$ has  infinite breadth, then there is a spread in $\Omega$, denoted by $\cE=(E_n)_{n\in{\mathbb N}}$, such that at least one of the following statements holds:
\begin{newnum}
\item $\{\ssx\cap\join(\cE)\colon \ssx\in \cS\} \supseteq\TMAX(\cE)$. 
\item $\{\ssx\cap\join(\cE)\colon \ssx\in \cS\}\supseteq\TMIN(\cE)$.
\item $\{\ssx\cap\join(\cE)\colon \ssx\in \cS\}\supseteq\TORT(\cE)$ and there is an $\cS$-decisive colouring of $\cE$.
\end{newnum}
\end{thm}
\begin{proof}
Let $\cS\subseteq\cP(\Omega)$ be a union-closed set system of infinite breadth.
It is enough to show that if we are in Case (S) of Theorem~\ref{t:dichotomy}, then there is a spread $\cE'$ in $\Omega$ such that $\left\{\ssx\cap\join(\cE') \colon \ssx\in \cS\right\}$ contains $\TMAX(\cE')$. This fact, together with Theorem~\ref{structure of infinite breadth} finishes the proof. 

Now assume that Case (S)  holds. So there exists a sequence $(\ssa_n)_{n\geq 1}$ of distinct members of $\cS$  so that for each $m\in\Nat$, the set $\{\ssa_1,\dots, \ssa_m\}$ is incompressible.
Fix positive integers $n_1<n_2<\dots$ such that $n_{k+1}-n_k \to\infty$ (for instance we could take $n_k=k^2$). Let $\ssd_k = \bigsetcup_{i=1}^{n_k} \ssa_i \in \cS$, and for  convenience set $\ssd_0=\emptyset$, $n_0=0$.
It is easy to see that for each $k\in\Nat$, the set
\[
\cF_k'\defeq \{ \ssa_j \setminus \ssd_{k-1} \colon n_{k-1}+1\leq j \leq n_k\}
\]
is an incompressible subset of the semilattice $\left\{\ssx\setminus\ssd_{k-1}: \ssx\in \cS\right\}$.
Moreover: since $\cF_k'$ is incompressible, for each $j$ with $n_{k-1}+1\leq j\leq n_k$, we can select an element of $\ssa_j\setminus \ssd_{k-1}$ that does not belong to any other member of $\cF_k'$. Let $E_k$ be the set of all these elements (very loosely, one can think of $E_k$ as a ``transversal'' for $\cF_k'$).
 Since $E_k\subseteq \ssd_k\setminus \ssd_{k-1}$ and $|E_k|= n_{k+1}-n_k$, the sequence $\cE=(E_k)_{k\in\Nat}$ is a spread in $\Omega$.

To finish, it suffices to show that given $k\in\Nat$ and some $\omega\in E_k$, there exists some $\ssz\in \cS$ such that $\ssz\cap \join(\cE)= E_{<k} \dotcup \{\omega\}$. By construction, there exists $\ssz'\in \cF'_k$ such that $\ssz'\cap E_k = \{\omega\}$, and also $E_1\cup\dots \cup E_{k-1}=\ssd_{k-1}\cap\join(\cE)$.
Put $\ssz=\ssz'\cup \ssd_{k-1}$: this satisfies $\ssz\cap\join(\cE) = E_1\cup\dots\cup E_{k-1}\cup\{\omega\}$, and we must show $\ssz\in\cS$. But since $\ssz'$ has the form $\ssa_i\setminus \ssd_{k-1}$ for some $i\in\Nat$, we have $\ssz=\ssa_i\cup \ssd_{k-1}\in\cS$, as required.
\end{proof}

\end{section}

\begin{section}{Constructing non-AMNM weights in the $\TORT$ case}
\label{sec:weight-case 3}
In this final case, by Theorem~\ref{structure of infinite breadth-refined}, we suppose that  there exists a spread $\cE$ with a colouring of $\cE$ which
 is $\cS$-decisive such that $\{x\cap\join(\cE): x\in \cS\} $ contains $\TORT(\cE)$.

\begin{lem}[Creating a log-weight from a decisive colouring]
\label{l:make-a-weight}
Let $\cS\subseteq\cP(\Omega)$ be union-closed, and let $\cE=(E_n)_{n\geq 1}$ be a spread in $\Omega$. Suppose there is an $\cS$-decisive colouring of $\cE$, call it $\cC$, with a decisive colour class $C_0$. 

For each $\ssx\in \cS$ define
\[ T(\ssx) = \{ n \in \Nat \colon |\ssx\setcap C \setcap E_n | \leq \frac{1}{2} |C\setcap E_n| \;\text{for all $C\in \cC$}\} \]
and define $\lm(\ssx) =\sup_{n\in T(\ssx)} | \ssx\setcap C_0 \setcap E_n |$. Then $\lm(\ssx)<\infty$, and $\lm : \cS \to [0,\infty)$ is a log-weight.
\end{lem}

Note that if $T(\ssx)=\emptyset$, then $\lm(\ssx)=0$, i.e.\ we take the usual convention when considering the least upper bounds of subsets of $[0,\infty)$. 

\begin{proof}
The first step is to show $\lm(\ssx) <\infty$. If $T(\ssx)$ is finite there is nothing to prove; so assume $T(\ssx)$ is infinite. Since $\cC$ is a colouring of the spread $\cE$, we have $\min_{C\in \cC} |C\setcap E_n | \to\infty$. Therefore, since $T(\ssx)$ is infinite,  $|\ssx^c\setcap C \setcap E_n | \to\infty$ along $n\in T(\ssx)$. Now it follows from the condition~\eqref{eq:decisive} and the definition of $\lm$ that $\lm(\ssx)<\infty$.

Finally, given $\ssx$ and $\ssy$ in $\cS$, observe that $T(\ssx\setcup \ssy) \subseteq T(\ssx)\setcap T(\ssy)$. Hence
\[ \lm(\ssx\setcup\ssy) \leq \sup_{n \in T(\ssx\setcup\ssy)} \left(|\ssx\setcap C_0\setcap E_n| + |\ssy\setcap C_0 \setcap E_n|\right) \leq \lm(\ssx)+\lm(\ssy),\]
as required.
\end{proof}

\begin{prop}
Let $\cS$, $\cE$, $\cC$ and $\lm$ be as in the previous lemma, and assume that $\left\{\ssx\cap\join(\cE) \colon \ssx\in \cS\right\}\supseteq \TORT(\cE)$. Then $(\cS,\lm)$ does not have $1$-propagation.
\end{prop}

\begin{proof}
Let $C_0$ be the decisive colour class used to define~$\lm$.
We fix $n\in\Nat$ and enumerate the elements of $E_n\setcap C_0$ as $\gamma_1,\dots,\gamma_M$,
where $n$ is sufficiently large so that $M\geq 2$.
Since 
$\cS$ restricted to $\join(\cE)$ contains $\TORT(\cE)$, there exist $\ssx_1,\dots,\ssx_M\in\cS$ such that
\[ \ssx_j \setcap \join(\cE)  =  E_{<n} \dotcup \{\gamma_j\}\dotcup E_{>n} \qquad(1\leq j\leq M). \]
In particular $\ssx_i\setcap E_n = \ssx_i\setcap C_0\setcap E_n = \{\gamma_i\}$, while $\ssx_i\setcap E_j=E_j$ for all $j\neq n$. It follows that $\lm(\ssx_i)=1$ for $i=1,\dots,M$.

Let $\cF_n=\{\ssx_1,\dots, \ssx_M\}\subseteq W_1$\/, and let $\ssb_n=\join(\cF_n)$. Since $\ssb_n\setcap E_n = C_0 \setcap E_n$, and since $\ssb_n\setcap E_m=E_m$ for all $m\neq n$, we have $\lm(\ssb_n)=0$. Hence $\ssb_n\in \filgen{\cF_n}\cap W_1$\/.

Let $K\geq 1$ and suppose $\ssb_n\in \bigcup_{m=0}^\infty\FBP_K^m(\cF_n)$. Then, in particular, there exists $m\geq 1$ with the following property:
\begin{equation}\label{eq:a property}
\text{there exists $\ssy\in\FBP_K^m(\cF_n)$ such that $|\ssy\setcap C_0\setcap E_n| > \frac{1}{2}|C_0\setcap E_n|$.}
\tag{$*$}
\end{equation}
Let $m$ be minimal with respect to this property, and let $\ssy$ satisfy~\eqref{eq:a property}. Then there exist $\ssy_1$ and $\ssy_2$ in $\FBP_K^{m-1}(\cF_n)$  such that $\ssy\subseteq \ssy_1\setcup\ssy_2$; again, when $m=1$, our convention here is that $\FBP_K^0(\cF_n)=\cF_n$. By the minimality of $m$ (when $m\ge 2$) or since $\abs{C_0\cap E_n}\ge 2$ and $\abs{\ssx_j\cap C_0\cap E_n}=1$ (when $m=1$),
\[
 	|\ssy_1\setcap C_0\setcap E_n| \le \frac{1}{2}|C_0\setcap E_n|\quad\text{and}\quad 
 |\ssy_2\setcap C_0\setcap E_n| \le \frac{1}{2}|C_0\setcap E_n|\;.
\]

Clearly, $\filgen{\cF_n}\subseteq \cP(\join(\cF_n))=\cP(\ssb_n)$. In particular, for $i\in\{1,2\}$, we have $\ssy_i\subseteq \ssb_n$ and hence
\[ \ssy_i\setcap E_n \subseteq \ssb_n\setcap E_n = C_0\setcap E_n.\]
This implies that $n\in T(\ssy_i)$, so that
\[ \lm(\ssy_i) = \sup_{j\in T(\ssy_i)} |\ssy_i\setcap C_0\setcap E_j| \geq |\ssy_i\setcap C_0\setcap E_n| \;. \]
Putting this all together, and remembering that $\ssy_1$ and $\ssy_2$ belong to $W_K$,
\[ \begin{aligned}
\frac{1}{2}|C_0\setcap E_n|
& \leq |(\ssy_1\setcup\ssy_2)\setcap C_0\setcap E_n| \\
& \leq |\ssy_1\setcap C_0\setcap E_n| +  |\ssy_2\setcap C_0\setcap E_n|\leq  \lm(\ssy_1)+\lm(\ssy_2)\leq 2K\,.
\end{aligned} \]
Hence $K \geq \frac{1}{4} |C_0\setcap E_n|$. It follows that $V_{\cF_n}(\ssb_n) \geq \frac{1}{4} |C_0\setcap E_n|\to \infty$, as required.
\end{proof}

\end{section}

\begin{section}{Final examples}
Our construction of non-AMNM weights could have been streamlined if the following claims were true:
\begin{newnum}
\item if $q:S\to T$ is a surjective homomorphism of semilattices, and $(T,\lm^*)$ fails $L$-propagation, then $(S,\lm^*\circ q)$ fails $L$-propagation;
\item if $(R,\lm)$ has $L$-propagation then so does $(T,\lm)$ for every subsemilattice $T\subset R$.
\end{newnum}
In this section we present examples to show that both claims are false.

Let $\Omega =\{0,1,2\}\times\Nat$. Consider two sequences $(\ssx_j)_{j=1}^\infty$ and $(\ssa_n)_{n=1}^\infty$ defined as follows:
\begin{equation}
\ssx_j \defeq \set{ (1,j) , (2,j) } \quad;\quad \ssa_n \defeq \set{ (0, k), (1,k) \colon n^2\leq k < (n+1)^2 }.
\end{equation}
It may be helpful to picture these sets as certain rectangular ``tiles'', as in Figure~\ref{fig2}.

\begin{figure}[hpt]
\hfil\includegraphics[scale=0.6]{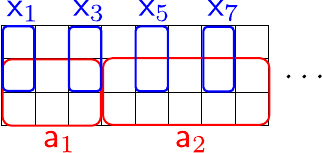}\hfil
\caption{Sets $\ssx_j$ and $\ssa_n$}
\label{fig2}
\end{figure}

Let $\cA=\set{ \ssa_n \colon n\in\Nat }$ and $\cX=\set{ \ssx_j \colon j\in \Nat }$. Define $\cS$ to be the semilattice generated inside $\cP(\Omega)$ by $\cA$ and $\cX$. A little thought shows that every member of $\cS$ has a unique decomposition, up to ordering, as a union of members of $\cA$ and members of $\cX$. We will think of the members of $\cA$ and $\cX$ as ``prime factors''.
Recall that for union-closed set systems, $\ssa$ is a factor of $\ssb$ precisely when $\ssa\subseteq\ssb$.

Now let $\Omega^* = \{1,2\}\times\Nat \subset \Omega$, and consider the \dt{truncation homomorphism}
$q: \cP(\Omega)\to\cP(\Omega^*)$ defined by $\ssz \mapsto \ssz\cap \Omega^*$.
Any log-weight $\lm^*$ defined on $\cP(\Omega^*)$ pulls back to give a log-weight $\lm=\lm^*\circ q$ on $\cP(\Omega)$. We now put
\begin{equation}
\lm^*(\ssz) \defeq |\{ j\in\Nat \colon (2,j) \in \ssz\}|
\end{equation}
which is clearly subadditive, and so gives us a log-weight on $\cP(\Omega^*)$. Note that $\lm$ is given by the same formula as $\lm^*$. Moreover, if $\ssz\in \cS$, then $\lm(\ssz)$ counts how many factors from $\cX$ occur in the ``prime factorization'' of $\ssz$.

\begin{eg}
Define $\cT=q(\cS)\subset \cP(\Omega^*)$. We will show that $(\cT,\lm^*)$ does not have $1$-propagation.

Let $\ssb_n\defeq q(\ssa_n) = \set{ (1,k) \colon n^2\leq k < (n+1)^2 }$ and let $\cB=\set{ \ssb_n \colon n\in\Nat}$. Then $\cT$ is the subsemilattice of $\cP(\Omega^*)$ generated by $\cX$ and $\cB$. For each $n$, let $\cE_n=\{\ssx_k \colon n^2\leq k< (n+1)^2\}$. Then $\cE_n\subseteq  W_1(\cT)$ and $\ssb_n \in \filgen{\cE_n}\cap W_1(\cT)$. Therefore, it suffices to prove that  $V_{\cE_n}(\ssb_n) \to\infty$ as $n\to\infty$.

\begin{lem}\label{l:another lemma}
Let $1\leq C \leq n$, and let $m\geq 0$. Then $\FBP_C^m(\cE_n)$ is contained in the union-closed set system generated by $\cE_n$, which we denote by $\langle \cE_n\rangle$.
\end{lem}

\begin{proof}
We induct on $m$.
The case $m=0$ is trivial. If the claim holds for $m=k-1$ where $k\in\Nat$, then let $\ssy, \ssz\in \FBP_C^{k-1}(\cE_n)$ and let $\ssa\in W_C(\cT)$ satisfy $\ssa\subseteq \ssy\setcup\ssz$. To complete the inductive step, we will show that $\ssa\in\langle\cE_n\rangle$.

By the inductive hypothesis, $\ssy\in\langle\cE_n\rangle$ and $\ssz\in\langle\cE_n\rangle$.
Therefore, by the ``prime factorization'' property of $\cS$, the only possible factors of $\ssa$ in $\cT$ are the members of $\cE_n$ together with $\ssb_n$. Suppose $\ssb_n\subseteq \ssa$; then $(1,k) \in \ssy\cup\ssz$ for every $n^2\leq k \leq n^2+2n$.
But this is impossible since $\lambda^*(\ssy)\leq C \leq n$ and $\lambda^*(\ssz)\leq C\leq n$.
\end{proof}

Since $\ssb_n\notin \langle\cE_n\rangle$, Lemma \ref{l:another lemma} implies that $V_{\cE_n}(\ssb_n) \geq n$, as required.
\end{eg}

\begin{lem}
$(\cS,\lm)$ has $L$-propagation for all $L\geq 0$.

\begin{proof}
Fix $L\geq 0$ and let $\cE$ be a non-empty subset of $W_L$.
Let $\fac(\cE)$ denote the set of factors of $\cE$ inside $\cS$.
Given $\ssz\in \filgen{\cE}\cap W_L$, we will show that $V_{\cE}(\ssz)\leq L$.

First note that each ``prime factor'' of $\ssz$ must be a factor of some element in $\cE$ (by unique factorization). Hence, $\ssz=\ssa'\cup\ssx'$ where $\ssa'$ is a product of members of $\cA\cap \fac(\cE)$ and $\ssx'$ is a product of members of $\cX\cap\fac(\cE)$.

Write $\ssx' = \ssx_{n(1)} \cup \ssx_{n(2)}\cup\dots \cup \ssx_{n(k)}$ where $n(1)<n(2)<\dots < n(k)$. Then $k=\lm(\ssz)\leq L$.  If $L<1$ then $\ssz=\ssa'$; induction on the number $m$ of the ``prime factors'' of $\ssa'$ yields $\ssa'\in \FBP^m_0(\cE)\subseteq \FBP^\infty_L(\cE)$, and so we are done. If $L\geq 1$, then by inductively considering $\ssy_0\defeq \ssa'$, $\ssy_1\defeq \ssy_0\cup\ssx_{n(1)}$, $\ssy_2 \defeq \ssy_1\cup\ssx_{n(2)}$, etc., we obtain $\ssy_j\in \FBP_L^{j+m}(\cE)$ for all $j=0,1,\dots, k$. Since $\ssz=\ssy_k\in\FBP_L^\infty(\cE)$, this completes the proof.
\end{proof}
\end{lem}

For our next example, consider the sets $\ssg_j \defeq \{( 1,j)\}$ for $j\in\Nat$. Define $\cG=\{ \ssg_j \colon j\in\Nat\}$ and define $\cR$ to be the semilattice generated by $\cX$ and $\cG$. Since each $\ssb_n$ is the union of finitely many members of $\cG$, $\cT$ is a subsemilattice of $\cR$.

\begin{eg}
One can check that $(\cR,\lm^*)$ has $L$-propagation for all $L\geq 0$. The proof is very similar to the proof for $(\cS,\lm)$. Unlike $\cS$, the semilattice $\cR$ does not have ``unique factorization''; but each $\ssz\in \cR$ has a largest factor belonging to $\langle\cX\rangle$, and this factor has a unique decomposition as a union of members of $\cX$. Using this, one can carry out the same kind of argument that was used to show $(\cS,\lm)$ has $L$-propagation. We leave the details to the reader.
\end{eg}

\end{section}


\section*{Acknowledgements}
This paper is the conclusion of a larger project, which grew out of conversations between the authors while attending the conference ``Banach Algebras and Applications'', held in Gothenburg, Sweden, July--August 2013, and was further developed while the authors were attending the thematic program ``Abstract Harmonic Analysis, Banach and Operator Algebras''  at the Fields Institute, Canada, during March--April 2014. The authors thank the organizers of these meetings for invitations to attend and for pleasant environments to discuss research.

The first author acknowledges the financial support of the Faculty of Science and Technology at Lancaster University, in the form of a travel grant to attend the latter meeting. The third author acknowledges the financial supports of a Fast Start Marsden Grant and of Victoria University of Wellington to attend these meetings.

The writing of this paper was facilitated by a visit of the first author to the University of Delaware in October 2022, supported by a Scheme 4 grant from the London Mathematical Society (reference 42128).
The second author also acknowledges support from National Science Foundation
Grant DMS-1902301 during the preparation of this article.


\newcommand{\etalchar}[1]{$^{#1}$}

\vfill

\newcommand{\address}[1]{{\small\sc#1.}}
\newcommand{\email}[1]{\texttt{#1}}

\noindent
\address{Yemon Choi,
School of Mathematical Sciences,
Lancaster University,
Lancaster LA1 4YF, United Kingdom} 

\email{y.choi1@lancaster.ac.uk}

\noindent
\address{Mahya Ghandehari,
Department of Mathematical Sciences,
University of Delaware,
Newark, Delaware 19716, United States of America}

\email{mahya@udel.edu}

\noindent
\address{Hung Le Pham,
School of Mathematics and Statistics,
Victoria University of Wellington,
Wellington 6140, New Zealand}

\email{hung.pham@vuw.ac.nz}

\end{document}